\documentclass[12pt,draft,amsmath,amscd,amsbsy,amssymb]{amsart}
 \relax

\chardef\bslash=`\\ 

\makeatletter \def\verbatim{\interlinepenalty\@M \@verbatim
\leftskip\@totalleftmargin\advance\leftskip2pc \frenchspacing\@vobeyspaces
\@xverbatim} \makeatother \hfuzz1pc

\makeatletter \def\dgt@k{\dg@DX=-3 \dg@DY=2 \dg@SIZE=3} \makeatother

\makeatletter \def\dgt@kk{\dg@DX=3 \dg@DY=-1 \dg@SIZE=3}

\theoremstyle{plain}
\newtheorem{thm}{Theorem}[section]

\newtheorem{lemma}[thm]{Lemma}
\newtheorem{prop}[thm]{Proposition}

\theoremstyle{definition}

\newtheorem{que}[thm]{Question}

\newcommand{\R}{{\mathbb R}}

\newcommand{\M}{{\mathcal M}}

\newcommand{\PM}{{\mathcal P}{\mathcal M}}

\numberwithin{equation}{section}

\newcommand{\diam}{{\mathrm d}{\mathrm i} {\mathrm a}{\mathrm m}\ }

\newcommand{\dom}{{\mathrm d}{\mathrm o}{\mathrm m}\ }

\newcommand{\N}{{\mathbb N}}

\newcommand{\A}{{\mathcal A}}

\newcommand{\e}{{\varepsilon}}

\begin{document}


\title[Linear extensions of partial (pseudo)metrics]
{On simultaneous linear extensions of partial (pseudo)metrics}
\author{E.D. Tymchatyn}

\address{Department of Mathematics and Statistics,
University of Saskatchewan,  Room 142 McLean Hall, 106 Wiggins Road,
  Saskatoon, SK. S7N 5E6 CANADA } \email{tymchat@math.usask.ca}

\author{M.~Zarichnyi}

\address{Department of Mechanics and Mathematics, Lviv National University,
Universytetska 1, 79000 Lviv, Ukraine} \email{topos@franko.lviv.ua,
mzar@litech.lviv.ua}
\thanks{The paper was finished during the second named author's visit to the
University of Saskatche\-wan.
This research was supported in part by NSERC research grant OGP 005616.}

\begin{abstract}    We consider the question of simultaneous extension of
(pseudo)me\-trics defined on nonempty closed subsets of a compact metrizable 
space. The main result is a counterpart of the  result due to 
K\"unzi and 
Shapiro for the case of  extension operators of partial continuous functions 
and
includes, as a special case, Banakh's theorem on linear regular operators 
extending
(pseudo)metrics.
\end{abstract}

\maketitle

\section{Introduction}\label{s:intro}

The theories of continuous extensions of continuous functions and continuous 
metrics develop in parallel.
Hausdorff's theorem on extension of a metric defined on a closed subset of a metrizable
space \cite{h} is a counterpart of the Tietze-Urysohn theorem on extensions of continuous functions.
Hausdorff's theorem has had many improvements; see for example 
\cite{bing, L1, t}.

The set of continuous pseudometrics on a compact metrizable $X$ forms a 
positive cone in the normed space $C(X\times X)$ of continuous functions on $X\times X$. 
C. Bessaga in \cite{be1, be2} formulated the problem  of existence of 
continuous linear operators 
that extend (pseudo)metrics defined on a closed subset of $X$. These
operators 
resemble the linear extension operators for continuous functions (see \cite{d}).
 In some special cases,
Bessaga solved this problem; his solution is based on an explicit formula
for extension of metrics onto the so called squeezed cone over a space. 
In its general setting, the 
extension problem for pseudometrics was first solved by T. Banakh \cite{b} by improving Bessaga's 
formulas and the proof was completed by O. Pikhurko \cite{p}. 
In \cite{b1}, T. Banakh  also developed another approach 
to the extension problem.  A very short 
proof of the existence of the linear operators extending (pseudo)metrics is
given in \cite{z}. Note that some sublinear extension operators for (pseudo)metrics
were constructed in \cite{nn}.

Recently, K\"unzi and Shapiro \cite{ks} considered a problem of simultaneous extension 
of partial continuous functions.  
The spaces of partial functions (i.e. the functions defined
on various subsets of topological spaces), first defined in 
\cite{k, k1}, naturally appear in the topological theory of 
differential equations (see \cite{ff}). In the present paper, we consider 
the analogous problem for partial 
(pseudo)metrics, i.e.  (pseudo)metrics defined on  closed subsets of a compact metrizable space.
Remark that the metric spaces with variable domains are of interest in different areas of 
geometry and topology. In particular, the (classes of equivalence up to 
isometry of) compact metric spaces are the points of the Hausdorff 
moduli  space (see, e.g. \cite{g}).

The main result of the paper is a theorem whose formulation is a mixture of those of
the cited theorems proved by Banakh and K\"unzi and Shapiro: There exists a  
continuous
operator that extends (pseudo)metrics defined on all closed subsets of a compact metrizable space
over all the space and whose restriction on the set of pseudometrics defined on 
every single subset is linear and regular (i.e. of norm 1; see the definitions below).

\section{Preliminaries}\label{s:prelim}

Let $X$ be a metrizable separable space. To avoid trivialities we always
assume that $|X|\ge2$. Let  $\Delta_X$  denote the diagonal of $X$.
Let $I=[0,1]$.

By $\exp X$ we denote the
{\em hyperspace} of $X$, i. e. the set of all nonempty compact subsets of $X$
endowed with the Vietoris topology. A base of this topology consists
of the sets of the form $$\langle U_1,\dots,U_k\rangle=\{A\in \exp X\mid
A\subset\cup_{i=1}^k U_i,\ A\cap U_i\neq\emptyset\text{ for all }i \},$$
where $ U_1,\dots,U_k$ run over the family of open subsets in $X$.

If $d$ is a compatible metric on $X$, then the Vietoris topology is
generated by the Hausdorff metric $d_H$,
$$d_H(A,B)=\inf\{\e>0\mid A\subset O_\e(B),\ B\subset  O_\e(A)\}.$$

Given a nonempty compact subset $A$ of $X$,
we denote by $\M(A)$ (respectively $\PM(A)$) the set of continuous metrics
(respectively pseudometrics) on $A$ (the fact  that $\varrho\in\PM(A)$ 
will  be also expressed as $\dom\varrho=A$). Set $\M=\cup\{\M(A)\mid A\in\exp X\}$,
$\PM=\cup\{\PM(A)\mid A\in\exp X\}$.

Identifying every pseudometric $d\in \PM$ with its graph,  which is a
compact subset of $X\times X\times\R$, we consider the set $\PM$ as a subset
of $\exp(X\times X\times\R)$ and endow $\PM$ with the subspace topology.

Let $K=\{(\varrho,\lambda)\in \PM\times\PM\mid \dom\varrho=\dom\lambda\}$. 
The space $\PM$ is a {\em positive cone} in the following sense: the maps 
$(c,\varrho)\mapsto c\varrho$ and $(\varrho,\lambda)\mapsto \varrho+\lambda$ are
continuous as maps from, respectively, $\R_+\times\PM$ and $K$ into $\PM$.

A map $u\colon \PM\to \PM(X)$ is called {\em linear} if (i)
$u(c\varrho)=cu(\varrho)$ for every $c\in\R_+$ and $\varrho\in \PM$ and (ii) 
$u(\varrho+\lambda)=u(\varrho)+u(\lambda)$ for every $(\varrho,\lambda)\in K$.

The {\em norm} of $\varrho\in\PM$ is $\|\varrho\|=\sup\{\varrho(x,y)\mid 
x,y\in \dom\varrho\}$.

A map $u\colon \PM\to \PM(X)$ is called {\em regular} if $\|u(\varrho)\|= \|\varrho\|$
for every $\varrho\in\PM$.

A map $u\colon \PM\to \PM(X)$ is called an {\em extension operator} if 
$u(\varrho)|(\dom\varrho\times\dom\varrho)=\varrho$, for every $\varrho\in\PM$.

\section{A selection theorem for multivalued mappings}\label{s:st}

Let $\langle T,\A,\mu\rangle$ be a measure space. For a Banach space 
$B$ with norm $\|\cdot\|$, 
we denote by $L_1(T,B)$ the Banach space of functions from $T$ to $B$ integrable in 
the Bochner sense. The norm in  $L_1(T,B)$ is defined by the formula
$\|\alpha\|=\int_T\|\alpha(t)\|d\mu$.

A set $Z$ of measurable
mappings from $\langle T,\A,\mu\rangle$  into a topological space $X$
is said to be {\em decomposable} (see, e. g., \cite{hu}, \cite{r}) if for every $f,g\in Z$ and every
$A\in\A$ the mapping
$$h(t)=\begin{cases} f(t), & t\in A,\\
g(t),& t\notin A
\end{cases}$$
belongs to $Z$.

Note that every set  in  $L_1(T,B)$ of cardinality 1 is decomposable. Besides, for every
closed subset $A$ of $B$ the set 
$L_1(T,A)=\{\alpha\in L_1(T,B)\mid \alpha(T)\subset A\}$ is closed decomposable for every
closed subset $A$ of $B$.

We will need the following selection theorem due to Fryszkowski \cite{f}.
Recall that a multi-valued mapping $F \colon X \to Y$ is called 
{\em lower semi-continuous} provided that 
$\{x\in  X \mid F(x) \cap U \neq \emptyset\}$ is open in $X$
 whenever $U$ is open in $Y$.

\begin{thm}\label{t:fry} Let $\langle T,\A,\mu\rangle$ be a compact measure space with a
non-atomic measure $\mu$ on a $\sigma$-algebra $\A$ of subsets of $T$ and $B$
a separable Banach space. Then every lower semicontinuous mapping $F$ from a
metric compactum $X$ into $L_1(T,B)$ with closed decomposable values admits
a continuous single-valued selection.
\end{thm}

We will apply this theorem in the case $T=I=[0,1]$ and $\mu=m$
(the Lebesgue measure on $I$). 

\section{Construction}\label{s:const}

Assume that a compact metric space $X$ is embedded into a
Banach space  $(B,\|\cdot\|)$. Define a multivalued map $F\colon \exp X\times X
\to L_1(I,B)$ as follows
$$F(A,x)=\begin{cases} \{x\}=L_1(I,\{x\}), &\text{ if }x\in A,\\
L_1(I,A), &\text{ if }x\notin A.\end{cases}$$

\begin{prop} The map $F$ is lower semicontinuous.
\end{prop}
\begin{proof} We have to show that for every open subset $U$ of $L_1(I,B)$ the
set $$U^\sharp=\{(A,x)\in \exp X\times X\mid F(A,x)\cap U\neq\emptyset\}$$
 is open in
$\exp X\times X$. Let $(A,x)\in U^\sharp$. Consider two cases.

Case 1). $x\notin A$. Then $F(A,x)=L_1(I,A)$. There exist $\alpha\in
L_1(X,A)\cap U$ and $\e>0$ such that $O_\e(\alpha)\subset U$. There exists a neighborhood $V$
of $(A,x)$ in $\exp X\times X$ such that $d_H(A,A')<\e/2$ and $x'\notin A'$
for every $(A',x')\in V$.

There exists $\beta\in L_1(I,A)$ such that $\|\alpha-\beta\|<\e/2$ and 
$\beta$ takes a finite number of values (see, e.g. \cite{DS}). 
One  can  perturb the values of 
$\beta$  to obtain  a function  $\beta'\in L_1(I,A')$ with the property that
$\|\beta(t)-\beta'(t)\|<\e/2$ for every $t\in I$.

 Since
$d(\alpha(t),\beta'(t))<\e$, for all $t$, we have
$\|\alpha-\beta'\|<\e$ and therefore $\beta'\in U$.
Hence, $V\subset U^\sharp$.

Case 2). $x\in A$. Then $F(A,x)=\{x\}$ (as above, we identify $x\in X$ with
 the constant function taking $x$ as its value). Since $x\subset U$ and 
$U$ is open, there exists $\e>0$ such that $O_\e(x)\subset U$. Then there is a neighborhood $V$
of $(A,x)$ in $\exp X\times X$ such that $d_H(A,A')<\e$ and $d(x,x')
<\e$ for every $(A',x')\in V$. For $(A',x')\in V$, denote by $\alpha'$ the constant map
whose value is $x'$ whenever $x'\in A'$ and arbitrary $y\in A'$ such that
$d(x,y)<\e$  whenever $x'\notin A'$. Then $x'\in F(A',x')\cap U$.

\end{proof}

Since $F$ satisfies the hypotheses of Theorem \ref{t:fry}, there exists a continuous
selection $f\colon \exp X\times X\to L_1(I,B)$   of $F$. Define $u=u_f\colon \PM\to\PM(X)$ by the formula
$$u(\varrho)(x,y)=\int_0^1\varrho(f(\dom(\varrho),x)(t),
f(\dom(\varrho),y)(t))dt.$$

\begin{lemma}\label{l:1} Let $\varrho\in\M(X)$. For every $\e>0$ there exits
$\delta>0$
such that \break $\int_0^1\varrho(\alpha(t),\beta(t))dt<\e$ whenever $\alpha,\beta
\in L_1(I,B)$ with $\|\alpha-
\beta\|<\delta$.
\end{lemma}

\begin{proof} We may assume that $|X|\ge2$. Since $\varrho$ is uniformly
continuous on $X\times X$, there exists $\delta_1>0$ such that
$\varrho(x,y)<\e/2$ whenever $\|x-y\|<\delta_1$.  Let
$0<\delta<\min\{ \delta_1,\frac{\e\delta_1}{2\diam X}\}$ (here $\diam X$ is the diameter of $X$ 
with respect  to the metric on $X$ induced by the norm in  $B$), then for every
$\alpha,\beta\in L_1(I,X)$ with $\|\alpha-\beta\|<\delta$ we have $m(A)<
\frac{\e}{2\diam X}$, where $A=\{t\mid
\|\alpha(t)-\beta(t)\|\ge\delta_1\})$ (recall that $m$  denotes the
Lebesgue measure on $I$). Therefore,
\begin{align*}
&\int_0^1\varrho(\alpha(t),\beta(t))dt=\int_A\varrho(\alpha(t),\beta(t))dt +
\int_{I\setminus A}\varrho(\alpha(t),\beta(t))dt\\ <&\frac{\e}{2\diam X}\diam X+
\int_{I\setminus A}\frac{\e}{2}dt\le\e.
\end{align*}
\end{proof}
\begin{lemma}\label{l:uni} Let   $\varrho\in\PM(X)$. The map $q=q_\varrho\colon L_1(I,X)\times L_1(I,X)\to\R$
defined by the formula $q(\alpha,\beta)=\int_0^1\varrho(\alpha(t),\beta(t))dt$,
is uniformly continuous.
\end{lemma}
\begin{proof}      Let $\e>0$. By Lemma \ref{l:1}, one can find $\delta>0$ such that
for every  $\alpha,\beta\in L_1(I,X)$ we have
 $\int_0^1\varrho(\alpha(t),\beta(t))dt<\frac{\e}{2}$ whenever $\|\alpha-
\beta\|<\delta$. Let now $\alpha_1,\beta_1,\alpha_2,\beta_2\in L_1(I,X)$
be such that $\|\alpha_1-\alpha_2\|<\delta$, $\|\beta_1-\beta_2\|<\delta$.
Then
\begin{align*}&|q(\alpha_1,\beta_1)-q(\alpha_2,\beta_2)|=
\left|\int_0^1(\varrho(\alpha_1(t),\beta_1(t))-
\varrho(\alpha_2(t),\beta_2(t)))dt\right|\\
\le & \left|\int_0^1(\varrho(\alpha_1(t),\beta_1(t))-
\varrho(\alpha_1(t),\beta_2(t)))dt\right|+\left|\int_0^1(\varrho(\alpha_1(t),\beta_2(t))-
\varrho(\alpha_2(t),\beta_2(t)))dt\right|\\
\le &\frac{\e}{2}+\frac{\e}{2}=\e.
 \end{align*}
 \end{proof}

\begin{prop}\label{p:contp} Let $\varrho\in\PM$. The function $u(\varrho)$ is 
a continuous pseudometric on $X$.
\end{prop}
\begin{proof}
It is evident that $u(\varrho)$ is a pseudometric on $X$. 

Let $(x_i)$, $(y_i)$ be convergent sequences in $X$ with the limits $x$
and $y$ respectively. Given $\e>0$, by Lemma \ref{l:uni}
one can find $\delta>0$ such that $|q_\varrho(\alpha_1,\beta_1)-
q_\varrho(\alpha_2,\beta_2)|<\e$ whenever $\|\alpha_1-\alpha_2\|<\delta$ and
 $\|\beta_1-\beta_2\|<\delta$ for 
 $\alpha_1,\alpha_2, \beta_1,\beta_2\in L_1(I,\dom\varrho)$. There exists $N\in\N$ such that for every $n>N$
$$\|f(\dom\varrho, x_n)-f(\dom\varrho, x)\|<\delta,\
\|f(\dom\varrho, x_n)-f(\dom\varrho, x)\|<\delta.  $$
Then $|u(\varrho)(x_n,y_n)-u(\varrho)(x,y)|<\e$ for every $n>N$.
\end{proof}

\begin{prop}\label{p:cont}
The map $u$ is continuous.
\end{prop}
\begin{proof} Let $(\varrho_n)$ be a sequence in $\PM$ and $\varrho_n\to\varrho\in\PM$.
Then, obviously, $\dom\varrho_n\to\dom\varrho$. Denote by $\tilde \varrho$
an extension of $\varrho$ over $X$. Denote by
$q\colon L_1(I,X)\times L_1(I,X)\to\R$ the map defined by the formula
$$q(\alpha,\beta)=\int_0^1\tilde\varrho(\alpha(t),\beta(t))dt,\ \alpha,\beta\in L_1(I,X).$$
As we have shown in Lemma \ref{l:uni}, $q$ is a uniformly continuous pseudometric on $L_1(I,X)$.

Let $\tilde \varrho_n=\tilde \varrho|(\dom\varrho_n\times \dom\varrho_n)$.
Obviously,
$\tilde\varrho_n\to\varrho$ in $\PM$ (see, e.g. \cite{ks}).

Let $\e>0$. We are going to show that there exists $n\in\N$ such that for every
$n\in\N$ with $n>N$ we have $|u(\varrho)(x,y)- u(\tilde\varrho_n)(x,y))|<\e$.

Since  $q$ is uniformly continuous, there exists $\delta>0$ such that for every
$\alpha,\beta\in L_1(I,X)$ with $\|\alpha-\beta\|<\delta$ we have
$q(\alpha,\beta)<\frac{\e}{2}$.

There exists $N_1\in\N$  such that for every $n\in\N$ with $n>N_1$ we have
$$\max\{|\tilde\varrho(x,y)-\varrho_n(x,y)|\mid (x,y)\in\dom\tilde\varrho_n
\}<\frac{\e}{2}.$$ Moreover, there exists $N_2\in\N$ such that for every
$n\in\N$ with $n>N_2$ and for every $x\in X$ we have
$\|f(\dom\varrho,x)-f(\dom\varrho_n,x)\|<\delta$.

Let $N=\max\{N_1,N_2\}$, $n>N$, and $x,y\in X$. Then
\begin{align*}
&|u(\varrho)(x,y)-u(\varrho_n)(x,y)|\\ = &\left|
\int_0^1\varrho(f(\dom\varrho,x)(t), f(\dom\varrho,y)(t))dt\right.\\
&\ \ \ \ \ \ \ \ \ \ \ \ \ \ \ \ \ \ \ \ \ \ \ \ \ \ \ \ \ \ \ \ \ \ \ \ \ \ \ \ 
\ \ \ \    -
\left.\int_0^1\varrho_n(f(\dom\varrho_n,x)(t), f(\dom\varrho_n,y)(t))dt\right|\\
\le&\left|\int_0^1\left(\tilde\varrho(f(\dom\varrho,x)(t), f(\dom\varrho,y)(t))-
\tilde\varrho(f(\dom\varrho_n,x)(t), f(\dom\varrho_n,y)(t) \right)dt\right|\\
+& \left|\int_0^1\left(\tilde\varrho(f(\dom\varrho_n,x)(t),
f(\dom\varrho_n,y)(t))- \varrho_n(f(\dom\varrho_n,x)(t),
f(\dom\varrho_n,y)(t) \right)dt\right|\\
\le &\frac{\e}{2}+\frac{\e}{2}=\e.
\end{align*}
\end{proof}
The main result of this section of this section is the following theorem which we improve 
in the forthcoming section.
\begin{thm}\label{t:2} There exists a continuous regular linear extension operator $u\colon \PM\to\PM(X)$.
\end{thm}
\begin{proof} It follows from Propositions \ref{p:contp} and \ref{p:cont}
 that $u\colon \PM\to\PM(X)$ is well-defined
and continuous. Obviously, $u$ is an extension operator. Its regularity and 
linearity are
easy consequences of the definition.

\end{proof}
\section{Operators preserving metrics}

The operator $u$ constructed in the previous section, in general, does not
preserve metrics, i.e., in general,  $u(\M)\not\subset\M(X)$. Note that there is
no linear operator $\tilde u\colon \PM\to\PM(X)$ such that $\tilde u(\M)\subset \M(X)$.
Indeed, for  every $a\in X$ we have $\PM(\{a\})=\M(\{a\})=\{0\}$ and it easily
follows from the linearity that $\tilde u(0)=0$.

Let $\tilde\M=\cup\{\M(A)\mid A\in \exp X,\ |A|\ge2\}$.

\begin{thm} There exists   a regular extension operator $\tilde u\colon \PM\to\PM(X)$
such that $\tilde u(\tilde\M)\subset\M(X)$.
\end{thm}
\begin{proof}
We are going to modify $u$ from Theorem \ref{t:2} as follows.

Given $(x,y)\in X\times X\setminus\Delta_X$ (as usual, $\Delta_X$ denotes the diagonal in $X\times X$)
and $A\in \exp X$ with $|A|\ge2$,
we define a multivalued
map $F_{(x,y,A)}\colon \exp X\times X\to L_1(I,X)$ as follows. Choose points
$a_x,a_y\in A$ so that
$a_x\neq a_y$ and $a_x=x$ (respectively $a_y=y$) if and only if $x\in A$
(respectively $y\in A$). Then
$$ F_{(x,y,A)}(B,z)=\begin{cases} F(B,z) & \text{ whenever }(A,x)\neq(B,z)\neq(A,y),\\
\{a_x\} & \text{ whenever }(B,z)=(A,x),\\
\{a_y\} & \text{ whenever }(B,z)=(A,y)
\end{cases}
$$
(here $F$ is the multivalued map from Section \ref{s:const}).

Note that the multivalued map $F_{(x,y,A)}$ is lower semicontinuous. Indeed, let $U$ be an open
subset of $L_1(I,X)$. Since $F$ is lower semicontinuous, the set $$U^\sharp=
\{(A',x')\in\exp X\times X\mid F(A',x')\cap U\neq\emptyset\}$$ is open in $\exp X\times X$.
Then the set
\begin{align*}&\{(A',x')\in\exp X\times X\mid F_{(x,y,A)}(A',x')\cap U\neq\emptyset\}\\
=&U^\sharp\setminus \{(A,z)\mid z\in\{x,y\},\ \{a_z\}\notin U\}
\end{align*}
is open in $\exp X\times X$.

Obviously, the images of  $F_{(x,y,A)}$ are closed and decomposable. By the Fryszkowski theorem, there exists a continuous selection $f_{(x,y,A)}$
of $F_{(x,y,A)}$. Since $f_{(x,y,A)}(A,x)\neq f_{(x,y,A)}(A,y)$, there exist
neighborhoods $W_A$ of $A$ in $\exp X$ and $V_{(x,y)}$ of $(x,y)$ in $X\times X$
respectively such that for every $A'\in W_A$, $(x',y')\in V_{(x,y)}$ we have
$f_{(x,y,A)}(A',x')\neq f_{(x,y,A)}(A',y')$.

Let $$Z=\{(x,y,A)\mid A\in\exp X,\ |A|\ge2,\ (x,y)\in
(X\times X)\setminus \Delta_X \}.$$   The cover
$\mathcal W=\{V_{(x,y)}\times W_A \mid (x,y,A)\in Z\}$ of $Z$ contains a countable
subcover $\mathcal W'= \{V_{(x_i,y_i)}\times W_{A_i} \mid i\in\N\}$ of $Z$.

For $i\in\N$ let $u_i=u_{f_{(x_i,y_i,A_i)}}$. Let $\tilde u=\sum_{i\in\N}\frac{u_i}{2^i}$.
Obviously, $\tilde u$ is a well-defined map from $\PM$ to $\PM(X)$. We leave to the reader 
a simple verification of the fact that $\tilde u$ is a regular extension operator of partial
pseudometrics on $X$.

We are going to show that $\tilde u(\tilde\M)\subset\M(X)$. Let $\varrho\in \tilde \M$
and $x,y\in X$, $x\neq y$. There exists $i\in\N$ such that
$(x,y,\dom\varrho)\in V_{(x_i,y_i)}\times W_{A_i}$. Then
$$u_i(\varrho)(x,y)=\int_0^1\varrho(f_{(x_i,y_i,A_i)}(\dom\varrho,x)(t),
f_{(x_i,y_i,A_i)}(\dom\varrho,y)(t))dt>0$$
and therefore $\tilde u(\varrho)(x,y)>0$. This proves that  $\tilde u(\varrho)$ is
a metric on $X$.

\end{proof}

\section{Non-metrizable case}

In the case of a non-metrizable compact Hausdorff space $X$, there are
no continuous operators of extension of 
partial continuous functions (E. Stepanova \cite{s}).

One can ask whether
there exists an extension operator $u\colon\PM\to\PM(X)$, for
 a non-metrizable compact Hausdorff space $X$. 

\begin{thm}\label{t:3} For a compact Hausdorff space $X$ the following conditions
are equivalent:
\begin{enumerate}
\item[(1)] there exists a continuous extension operator  $u\colon\PM\to\PM(X)$;
\item[(2)] there exists a continuous map 
$\Psi\colon (X\times X)\setminus\Delta_X\to \PM(X)$,
 $(x,y)\mapsto \Psi_{(x,y)}$, with $\Psi_{(x,y)}(x,y)\neq 0$
 for all
$(x,y)\in  X^2\setminus\Delta_X$;
\item[(3)]  $X$ is metrizable.
\end{enumerate}
\end{thm}

\begin{proof}

$(1)\Rightarrow(2)$. Let $u\colon \PM\to\PM(X)$  be a continuous extension operator.
Denote by $(x,y)\mapsto v_{(x,y)}\colon X^2\setminus\Delta_X\to \PM$ 
the map defined 
by  the conditions $\dom v_{(x,y)}=\{x,y\}$, $v_{(x,y)}(x,y)=1$.
Then the map $\Psi \colon X^2\setminus\Delta_X\to \PM(X)$,
$(x,y)\mapsto \Psi_{(x,y)}=u(v_{(x,y)})$ has the property that
$\Psi_{(x,y)}(x,y)=1$.

$(2)\Rightarrow(3)$. Given a  continuous map 
$\Psi\colon (X\times X)\setminus\Delta_X\to \PM(X)$,
 $(x,y)\mapsto \Psi_{(x,y)}$, with $\Psi_{(x,y)}(x,y)\neq 0$
 for all
$(x,y)\in  X^2\setminus\Delta_X$ define 
$\Phi\colon (X\times X)\setminus\Delta_X\to C(X)$ by the formula 
$\Phi_{(x,y)}(z)=\Psi_{x,y}(x,z)$, 
$z\in X$ . Then $\Phi_{(x,y)}(x)=0\neq 1=\Phi_{(x,y)}(y)$ 
 for all
$(x,y)\in  X^2\setminus\Delta_X$.
Therefore
the map $\Phi$ is a continuous separating map in the sense of \cite{s}. By
\cite[Proposition 2]{s}, the space $X$ is metrizable.

$(3)\Rightarrow(1)$. Follows from Theorem \ref{t:2}. 

\end{proof}

\section{Remarks and open questions}

The set $\mathbf E(X)$ of all regular extension operators of partial (pseudo)metrics on $X$ is a
topological invariant of $X$. One can  consider this set 
in  different topologies, say, in the uniform convergence topology, compact-open topology,
and pointwise convergence topology. 
\begin{que} If $\mathbf E(X_1)$ and $\mathbf E(X_2)$ are (topologically) isomorphic
with respect to the natural convex structure in them (and one of the topologies mentioned above),
are the compact metrizable spaces $X_1$ and $X_2$ homeomorphic?
\end{que}

The following question is motivated by the results of \cite{bb}.
\begin{que} Are there linear extension operators $u\colon\PM\to\PM(X)$ of norm 1 that 
are also
continuous in the  pointwise convergence topology on every $\PM(A)$, 
where $A\in\exp X$? 
\end{que}

Note that the implication $(2)\Rightarrow(3)$ from Theorem \ref{t:3} can be  
generalized on the class of paracompact $p$-spaces
(a topological space is a paracompact $p$-space if it admits a perfect map 
onto a metric space). The proof remains the same.

\end{document}